\documentclass[a4paper,12pt]{scrartcl}

\usepackage[bookmarks,pagebackref]{hyperref}
\usepackage{amsrefs}
\usepackage[utf8]{inputenc}
\usepackage[english]{babel}
\usepackage{amssymb}
\usepackage{amsmath}
\usepackage{latexsym}
\usepackage{amsthm}
\usepackage{eucal}
\usepackage{bbm}
\usepackage{chngcntr}
\usepackage{apptools}
\usepackage{mathtools}
\usepackage{authblk}
\usepackage[pdflatex]{crop}
\usepackage{tikz}
\usepackage{tikz-cd}
\usepackage{microtype}
\usepackage{enumitem}
\usepackage{graphicx}
\usepackage{mathrsfs}
\usepackage[colorinlistoftodos]{todonotes}
\usepackage{yfonts}
\usepackage{tabularx}
\usepackage{color}

\allowdisplaybreaks

\setkomafont{disposition}{\normalfont\bfseries}
\newcommand{\auth}[0]{{Renan Assimos and J\"urgen Jost}}
\newcommand{\tit}[0]{{The Geometry of  Maximum Principles and a Bernstein Theorem in Codimension 2}}
\newcommand{\kw}[0]{{Bernstein theorem, minimal graph, harmonic map, Grassmannian, Gauss map, maximum principle}}

\hypersetup{
pdfauthor={\auth},%
pdftitle={\tit},%
colorlinks, linktocpage=true, pdfstartpage=1, pdfstartview=FitV,%
breaklinks=true, pdfpagemode=UseNone, pageanchor=true, pdfpagemode=UseOutlines,%
plainpages=false, bookmarksnumbered, bookmarksopen=true, bookmarksopenlevel=1,%
hypertexnames=true, pdfhighlight=/O,%
urlcolor=black, linkcolor=black, citecolor=black, 
}
\numberwithin{equation}{section}

\theoremstyle{plain}
\newtheorem{thm}{Theorem}[section]
\AtAppendix{\counterwithin{thm}{section}}
\newtheorem{defn}[thm]{Definition}
\newtheorem{lemma}[thm]{Lemma}
\newtheorem{prop}[thm]{Proposition}

\newtheorem{notation}[thm]{Notation}

\newtheorem{condition}{Condition}
\newtheorem{rmk}[thm]{Remark}
\newtheorem{remark}[thm]{Remark}

\theoremstyle{definition}

\newtheorem{eg}[thm]{Example}

\DeclareMathOperator{\test}{span}

\newcommand{\Z}{\mathbb{Z}}
\newcommand{\N}{\mathbb{N}}

\newcommand{\R}{\mathbb{R}}

\DeclareMathOperator{\tr}{tr}

\DeclareMathOperator{\grafico}{graph}
\DeclareMathOperator{\Hem}{Hem}

\let\originalleft\left
\let\originalright\right
\renewcommand{\left}{\mathopen{}\mathclose\bgroup\originalleft}
\renewcommand{\right}{\aftergroup\egroup\originalright}

\newcommand{\ie}{\emph{i.e.} }

\title{\tit}
\author{\auth\thanks{Correspondence: \href{mailto:assimos@mis.mpg.de}{assimos@mis.mpg.de}, \href{mailto:jjost@mis.mpg.de}{jjost@mis.mpg.de}}}
\affil{\small Max Planck Institute for Mathematics in the Sciences\\ Leipzig, Germany}
\date{}

\begin{document}

\maketitle

\begin{abstract}
Moser's Bernstein theorem \cite{moser61} says that an entire minimal graph of codimension 1 with bounded slope must be a hyperplane. An analogous result for arbitrary codimension is not true, by an example of Lawson-Osserman. Here, we show that Moser's theorem nevertheless extends to codimension 2,  i.e., a  minimal $p$-submanifold in $R^{p+2}$, which is the graph of a smooth function defined on the entire $R^p$ with bounded slope, must be a $p$-plane. Our method depends on convexity properties of Grassmannians which come into play as the targets of  the (harmonic) Gauss maps of our minimal submanifolds. In fact, we develop a general method to construct subsets of complete Riemannian manifolds that cannot contain images of non-constant harmonic maps from compact manifolds. When applied to Grassmannians, for codimension 2, it yields an appropriate domain that cannot contain nontrivial images of such harmonic maps. Our conclusion will then be reached by an application of Allard's theorem to handle the issue that the minimal submanifolds in question are not compact. 
\end{abstract}

\textbf{Keywords: }{Bernstein theorem, minimal graph, harmonic map, Grassmannian, Gauss map, maximum principle.}

\newpage
\tableofcontents

\section{Introduction}

One of the central results  of the theory of minimal surfaces is Bernstein's theorem, stating that the only entire minimal graphs in Euclidean 3-space are planes. In other words, if $f(x,y)$ is a smooth function defined on all of $\R^2$ whose graph in $\R^3$, $(x,y,f(x,y))$, is a minimal surface, then $f$ is a linear function and its graph is a plane.

Profound methods in analysis and geometric measure theory were developed to generalize Bernstein's theorem to higher dimensions, culminating in the theorem of J. Simons \cite{simons68} stating that an entire minimal graph has to be planar for dimension $d\leq 7$. This dimension constraint is optimal, as  Bombieri, de Giorgi, and Giusti \cite{bombieri69} constructed  a counter-example to such an assertion in dimension 8 and higher. This reveals the subtlety and difficulty of the problem. Under the additional assumption that the slope of the graph is uniformly bounded, Moser \cite{moser61} proved a Bernstein-type result in arbitrary dimension. 

All the preceding results consider minimal hypersurfaces, that is, minimal graphs in Euclidean space of codimension 1. 
For higher codimension, the situation is more complicated. On one hand, Lawson-Osserman \cite{lawson77} have given explicit counterexamples to Bernstein-type results in higher codimension. Namely, the cone over a Hopf map is an entire Lipschitz solution to the minimal surface system. Since the slope of the graph of such a cone is bounded, even a Moser-type result for codimension higher than one cannot hold.

On the other hand, there are also some positive results in higher codimension, although, in view of the Lawson-Osserman examples, they necessarily require additional constraints. 
We can describe the main development as a sequence of steps. Those results all depend on  the fact that, by a theorem of Ruh-Vilms \cite{ruh70}, the Gauss map of a minimal submanifold in Euclidean space is harmonic. This Gauss map takes its values in a Grassmann manifold $G^{+}_{p,n}$ (which is a sphere in the case of codimension $n-p=1$). Therefore, the geometry of Grassmann manifolds is the key to understanding the scope of Bernstein theorems in higher codimension. Since the composition of a harmonic map (such as the Gauss map) with a convex function is a subharmonic function, when such a convex function is found the maximum principle can be applied to show that, when the domain of the harmonic map is compact, the resulting subharmonic function is constant. And when the convex function is nontrivial, for instance  strictly convex, then the harmonic map itself is constant, and the minimal graph is therefore linear. A key technical problem emerges from the fact that in our application, the domain is $\R^p$, which is not compact, so that the maximum principle cannot be applied directly. We postpone the discussion of this issue and return to the geometric steps. 
\begin{enumerate}
\item[1.] Hildebrandt-Jost-Widman \cite{hildebrandt80} identified the largest ball in the Grassmannian on which the squared distance function from the center is convex. Thus, when the Gauss image is contained in such a ball, that is, when the slopes of the tangent planes of our minimal submanifold do not deviate too much from a given direction, the minimal graph is linear, that is, a Bernstein result holds. 

\item[2.] Why consider only distance balls? Jost-Xin \cite{jost99} constructed regions in $G^{+}_{p,n}$ which are larger than convex metric balls but on which squared distance function is still convex. After all, even though $G^{+}_{p,n}$ is a symmetric space, the distance function behaves differently in different directions. Thus, the Gauss region implying the Bernstein results is larger.

\item[3.] Why consider only the squared distance function? Jost-Xin-Yang \cite{jost2012},\cite{jost2016} went further by constructing larger regions in $G^{+}_{p,n}$ that support strictly convex functions. Thus, the Bernstein theorem was further extended. In particular in codimension 1, that is, the classical case, even minimal hypersurfaces that might be  more general than graphs are shown to be linear. Still, it is not clear whether the Lawson-Osserman example is sharp or whether there exist other examples that come closer to the bounds obtained by \cite{jost2016} in higher codimension. 

\item[4.] The level sets of convex functions are convex. The starting idea of the present paper is that this is the key property: to have a family of convex hypersurfaces. But why do we need convex functions? We show that a foliation by convex hypersurfaces suffices, which does not need to come from a convex function. And this clarifies the geometric nature of the maximum principle that was the cornerstone of the reasoning just described. As far as we can tell, conceptually, this seems to be the ultimate step in this line of research. Theorem \ref{codim2}, to be stated shortly, seems to be the optimal result in codimension 2. It remains to explore our scheme in higher codimension. 
\end{enumerate}
In fact, all those results apply more generally to graphs of parallel mean curvature, since the Gauss map remains harmonic in such cases by the Ruh-Vilms theorem. It was proven by Chern \cite{chern65} that a hypersurface in Euclidean space that is an entire graph of constant mean curvature necessarily is a minimal hypersurface. Thus, by Simons' result, it is a hyperplane for dimension $d\leq 7$. See also Chen-Xin \cite{chen92} for a generalization of Chern's result.\\

A maximum principle, however, only applies directly when the domain is compact, but the domain of an entire minimal graph is $\R^p$. Therefore, one needs to turn the qualitative maximum principle into quantitative Harnack-type estimates, a technique also pioneered by Moser \cite{moser61}. In the proof of Hildebrandt-Jost-Widman \cite{hildebrandt80}, the analytical properties of such convex functions were used to derive  H\"older estimates for harmonic maps with values in Riemannian manifolds with an upper bound for the sectional curvature. By a scaling argument, they could conclude a Liouville-type theorem for harmonic maps under assumptions including the mentioned harmonic Gauss map. In the same setting, Jost-Xin-Yang refined the tools developed in \cite{hildebrandt80} and \cite{jost99} to obtain a-priori estimates for harmonic maps, improving higher codimension Bernstein results. 
Here, we use a geometric maximum principle due to Sampson \cite{sampson78} (see also \cite{eells95}), to develop a general method for constructing subsets  that do not admit images of non-constant harmonic maps defined on compact manifolds. This method does not allow us to obtain H\"older estimates, but fortunately, we can replace them by the result of Allard \cite{allard72}, which is a seminal results in geometric measure theory,   to study the graph case and  obtain Bernstein-type results. For our purposes, Allard's Theorem reduces the case of minimal submanifolds of Euclidean space to that of minimal submanifolds of spheres. The corresponding Gauss map for minimal submanifolds of spheres is still harmonic, so that the reasoning just described still works.

More concretely, while  Lawson-Osserman cones appear in codimension 3 or higher, we  prove a  Moser-type  result for codimension 2. 

\begin{thm}[Moser's Theorem in codim 2]\label{codim2}
	Let $z^i=f^i(x^1,...,x^p)$, $i=1,2$, be smooth functions defined everywhere in $\R^p$. Suppose the graph $M=(x,f(x))$ is a submanifold with parallel mean curvature in $\R^{p+2}$. Suppose that there exists a number $\beta_0<+\infty$ such that
    \begin{equation}
      \Delta_{f}(x)\leq \beta_0 \hspace{0.3cm}\text{for}\hspace{0.2cm}\text{all}\hspace{0.2cm} x\in\R^p,
    \end{equation}
    where
    \begin{equation}
      \Delta_{f}(x):=\big\{\det\big(\delta_{\alpha\beta} + \sum_{i}f^{i}_{x_{\alpha}}(x)f^{i}_{x_{\beta}}(x)\big)\big\}^{\frac{1}{2}}.
    \end{equation}
    Then $f^1,f^2$ are linear functions on $\R^p$ representing an affine p-plane in $\R^{p+2}$.
\end{thm}

The results of this paper constitute the main part of the first author's PhD thesis. We thank Lei Liu and Ruijun Wu for helpful discussions as well as Aleksander Klimek, Paolo Perrone and Sharwin Rezagholi for helpful comments. We also thank Luciano Mari for helping us to clarify some arguments on the previous version of this pre-print.

\section{Preliminaries}

Let $(M,g)$ and $(N,h)$ be Riemannian manifolds without boundaries. By Nash's Theorem we have an isometric embedding $N\hookrightarrow \R^L$. 

\begin{defn}
    A map $\phi \in W^{1,2}(M,N)$ is called harmonic iff it is a critical point of the energy functional\\
	\begin{equation}
	E(\phi) := \frac{1}{2} \int_{M} \Vert d\phi\Vert^{2}dvol_{g},
	\end{equation}
	where $\Vert .\Vert^{2} = \langle . ,.\rangle$ is the metric over the bundle $T^{*}M\otimes \phi^{-1}TN$ induced by $g$ and $h$.
\end{defn}

\noindent Recall that the Sobolev space $W^{1,2}(M,N)$ is defined as:
\begin{align}
W^{1,2}(M,N) :=\hspace{0.2cm} &\Big\{v: M\longrightarrow \R^L; \hspace{0.2cm}\left|\left|v\right|\right|^{2}_{W^{1,2}(M)}=\int_{M}(\left|v\right|^2 + \Vert dv\Vert^{2})\hspace{0.1cm}dv_g< +\infty \hspace{0.1cm}\text{and}\hspace{0.1cm}\\
 &v(x)\in N \hspace{0.1cm}\text{for}\hspace{0.1cm}\text{almost}\hspace{0.1cm}\text{every}\hspace{0.1cm} x\in M \Big\}.
\end{align}

The Euler-Lagrange equations for the energy functional are:
\begin{equation}
\tau(\phi)=0,
\end{equation}
and $\tau$ is called the tension field of the map $\phi$.

In local coordinates
\begin{align}
e(\phi)=\Vert d\phi\Vert^2 = g^{ij}\frac{\partial\phi^{\beta}}{\partial x^i}\frac{\partial\phi^{\gamma}}{\partial x^j}h_{\beta\gamma}.\\
\tau(\phi) = \big(\Delta_{g}\phi^{\alpha} + g^{ij}\Gamma_{\beta\gamma}^{\alpha}\frac{\partial\phi^{\beta}}{\partial x^i}\frac{\partial\phi^{\gamma}}{\partial x^j}h_{\beta\gamma}\big)\frac{\partial}{\partial\phi^{\alpha}}.
\end{align} 
where $\Gamma_{\beta\gamma}^{\alpha}$ denote the Christoffel symbols of N.

\begin{defn}[Gauss Map]
	Let $M^p\hookrightarrow \R^n$ be a p-dimensional oriented submanifold in Euclidean space. For any $x\in M$, by parallel translation in $\R^n$, the tangent space $T_{x}M$ can be moved to the origin, obtaining a p-subspace of $\R^n$, i.e., a point in the oriented Grassmannian manifold $G^{+}_{p,n}$. This defines a map $\gamma: M\longrightarrow G^{+}_{p,n}$ called the Gauss map of the embedding $M\hookrightarrow \R^n$.
\end{defn}

\begin{thm}[Ruh-Vilms]
	Let M be a submanifold in $R^n$ and let $\gamma:M\longrightarrow G^{+}_{p,n}$ be its Gauss map. Then $\gamma$ is harmonic if and only if M has parallel mean curvature.
\end{thm}

\noindent Harmonic maps have interesting geometric properties. By using Ruh-Vilms Theorem, one can try to find subsets $A\subset G^{+}_{p,n}$ for which there can be no non-constant harmonic map $\phi$ defined on some compact manifold M with $\phi(M)\subset A$. In that regard, it is often useful to use  the composition formula  for $\phi:M\longrightarrow N$, $\psi:N\longrightarrow P$ where $(P,i)$ is another Riemannian manifold, 
\begin{equation}
\tau(\psi\circ\phi) = d\psi\circ\tau(\phi) + \tr\nabla d\psi(d\phi,d\phi). 
\end{equation} 
When $\phi$ is harmonic, \ie $\tau(\phi)=0$, the formula is particularly useful. 
In particular if $P=\R$, and $\psi$ is a (strictly) convex function, then $\tau(\psi\circ\phi) \geq 0$ ($>0$). That is, $\psi\circ\phi:M\longrightarrow \R$ is a (strictly) subharmonic function on M. The maximum principle then implies the following proposition:

\begin{prop}
	Let $M$ be a compact manifold without boundary, $\phi:M\longrightarrow N$ a harmonic map with $\phi(M)\subset V \subset N$. Assume that there exists a strictly convex function on $V$. Then $\phi$ is a constant map.
\end{prop} 
\noindent In our setting, the target $N$ is the Grassmannian $G^{+}_{p,n}$. 
Then to obtain such a subset $A\subset G^{+}_{p,n}$, one tries to find a strictly convex function $f:A\longrightarrow \R$. This strategy was used by Hildebrandt-Jost-Widman, Jost-Xin, Jost-Xin-Yang and others \cite{hildebrandt80}, \cite{jost99}, \cite{jost2016}.

In this paper we want to modify this approach. Instead of using strong analytical arguments to obtain a subset that admits a strictly convex function, we want to explore the geometry of regions that can contain the image of a non-constant harmonic map.

To economize on notation, we state the following definition to be used throughout the paper.
 
\begin{defn}[Property ($\star$)]\label{star}
	 We say that an open connected subset $\mathcal{R}\subset (N,h)$, where $(N,h)$ is a complete Riemannian manifold, has {property ($\star$)} if there is no pair $(M,\phi)$, where $(M,g)$ is a compact manifold and $\phi:M \longrightarrow N$ is a non-constant harmonic map with $\phi(M)\subset \mathcal{R}$.
\end{defn}  

\section{The Geometry of the Maximum Principle}

\subsection{Sampson's Maximum Principle}

A beautiful result in the theory of harmonic maps is  Sampson's maximum principle (SMP):

\begin{thm}[SMP]\label{SMP}
 Let $\phi : (M,g) \longrightarrow (N,h)$ be a non-constant harmonic map, where $M$ is a compact Riemannian manifold, $N$ is a complete Riemannian manifold, and $S\subset N$ is a hypersurface with definite second fundamental form at a point $y=\phi(x)$. Then no neighborhood of $x\in M$ is mapped entirely to the concave side of $S$.
\end{thm}

\noindent A proof can be found in \cite{sampson78} and \cite{eells95}.

\begin{remark}
Take a geodesic ball $B(p,r)$ in a complete manifold $N$ such that $r$ is smaller than the convexity radius of $N$ at $p$. Then $\partial B(p,r)$ is a hypersurface of $N$ with definite second fundamental form for every point $q\in\partial B(p,r)$. 	
\end{remark}

\noindent The main result of this work can be seen as a corollary of Sampson's maximum principle.

\begin{thm}\label{corollary of SMP}
	Let $(N,h)$ be a complete Riemannian manifold and $\Gamma:[a,b]\longrightarrow N$ a smooth embedded curve. Consider a smooth function $r:[a,b]\longrightarrow \R_+$ and a region
	\begin{equation}\mathcal{R}:=\bigcup_{t\in[a,b]}B(\Gamma(t),r(t)),\end{equation}
	where $B(\cdot,\cdot)$ is the geodesic ball and $r(t)$ is smaller than the convexity radius of N for any $t$.
	If, for each $t_{0}\in (a,b)$, the set $\mathcal{R} \backslash B(\Gamma(t_{0}),r(t_{0}))$ is the union of two disjoint connected sets, namely the connected component of $\Gamma(a)$ and the one of $\Gamma(b)$, then there exists no compact manifold $(M,g)$ and non-constant harmonic map $\phi:M\longrightarrow N$ such that $\phi(M)\subset\mathcal{R}$. In other words, $\mathcal{R}$ has the property ($\star$) of Definition~\ref{star}.
\end{thm} 

\begin{figure}
	\begin{picture}(100,270)
	\put(0,0){\includegraphics[width=1\linewidth]{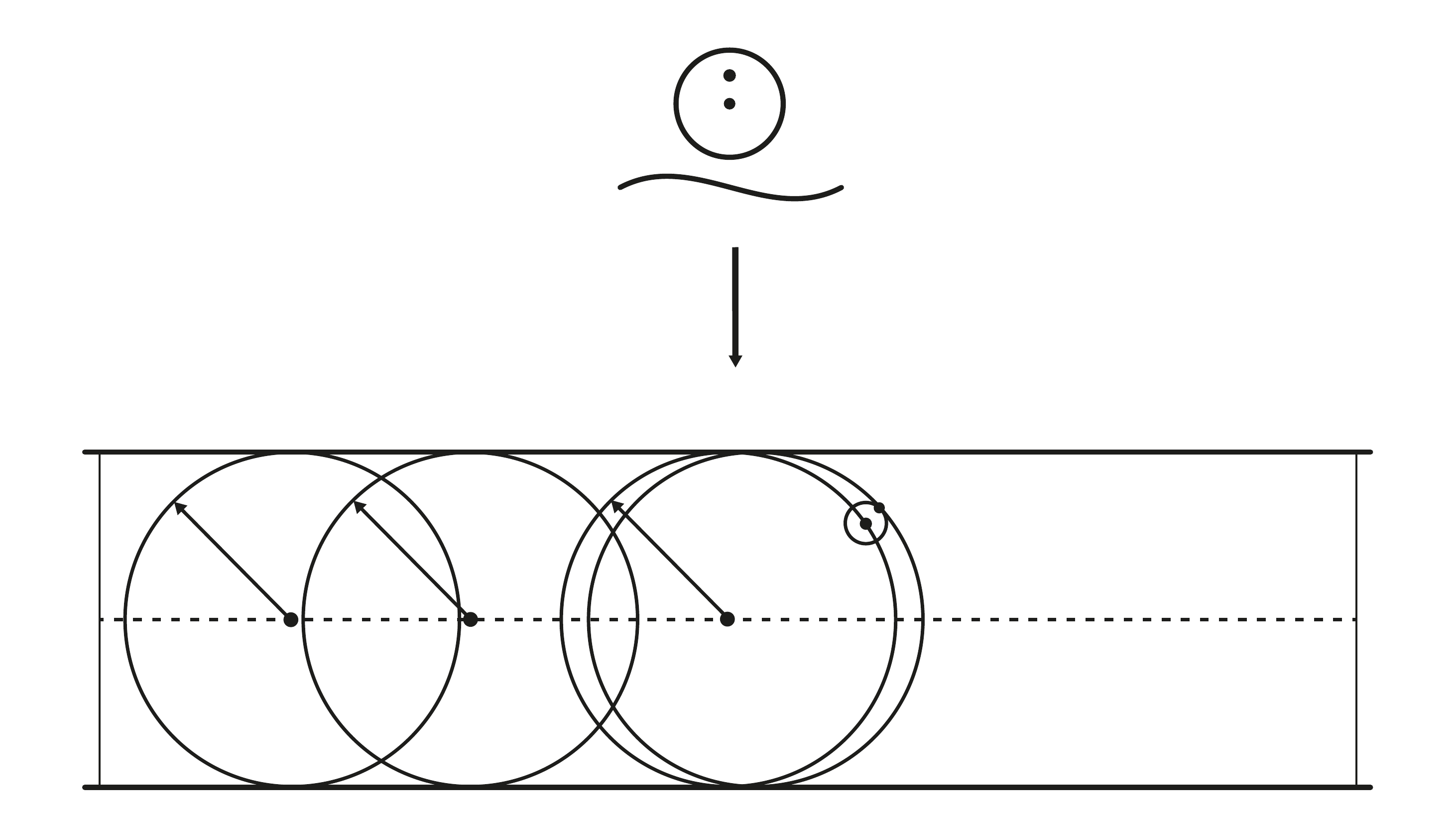}}
	\put(210,155){$\phi$}
	\put(235,155){harmonic}
	\put(280,205){$M$ cpt}
	\put(280,190){$\partial M=\emptyset$}
	\put(45,125){$\mathcal{R}:=\bigcup_{t\in[a,b]}B_g(\Gamma(t),r)$}
	\put(70,85){$r$}
	\put(125,85){$r$}
	\put(205,85){$r$}
	\put(275,97){$y_1=\phi(x_1)$}
	\put(240,78){$\phi(x)=y$}
	\put(227,225){$x_1$}
	\put(215,212){$x$}
	\put(5,57){$\Gamma(a)$}
	\put(340,65){$\Gamma(t)$}
	\put(425,55){$\Gamma(b)$}
	\end{picture}
	\caption{A region $\mathcal{R}$ as the union of convex balls}
	\label{Fig1}
\end{figure}
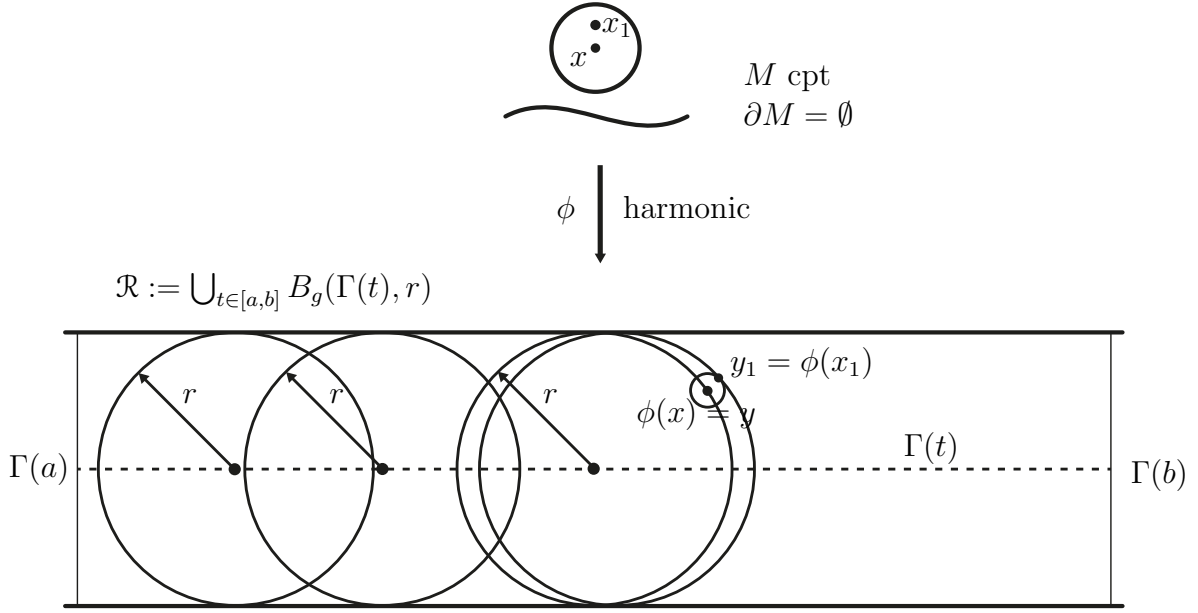

\noindent Before presenting the proof of this statement, we make two remarks on the geometry of the region $\mathcal{R}$

\begin{rmk}
	By hypothesis, for each time $t_0\in [a,b]$, the set $\mathcal{R}\backslash B(\Gamma(t_0),r(t_0))$ is the union of two disjoint connected components. This condition implicitly rules out a lot of possibilities for both $\Gamma$ and $r$. Namely, suppose we have a manifold with large convex radius. We could define $\mathcal{R}$ as the union of balls over a geodesic of unit speed $\Gamma$ of a very small length $\epsilon$ with radius varying smoothly from $\epsilon/2$ at $t=0$, to $10\epsilon$ at $t=\epsilon/2$, and back to $\epsilon/2$ when $t=\epsilon$. This region does not satisfies the necessary condition for our theorem.
	
	On the other hand, suppose we have a curve $\Gamma$ with the shape of a very steep `$U$'. Our condition implies that, in this case, the radius of the balls near the turning point of the `U' must be very small, otherwise there will be some intersection between the connected components of $\Gamma(a)$ and $\Gamma(b)$.
\end{rmk}

\begin{lemma}\label{property of R}
	Let $\mathcal{R}$ be a region defined as in Theorem~\ref{corollary of SMP} satifying the same hypothesis. Let $t_0\in (a,b)$ and  consider the open set
	\begin{equation}\label{connected component of gamma b}
	C_b(t_0):=\bigcup_{t\in\left(t_0,b\right)}B(\Gamma(t),r(t))\backslash \left(\overline{B(\Gamma(t_0),r(t_0))}\right)
	\end{equation}
	This is the connected component of $\Gamma(b)$ in $\mathcal{R}\backslash\left(B(\Gamma(t_0),r(t_0))\right)$. Then, for every point $y\in C_b(t_0)$, there exists $t_1>t_0$ such that $y\in\overline{C_b(t_1)}$ and $y\in \partial B(\Gamma(t_1),r(t_1))$. 
\end{lemma}

\begin{proof}
	By definition of $C_b(t_0)$, whenever $s'>s$, we have $C_b(s')\subsetneq C_b(s)$. Then for a fixed $y\in C_b(t_0)$, we define the smooth function $\kappa_y:[a,b]\longrightarrow \R$ given by $\kappa_y(s):=d\left(\Gamma(s),y\right)-r(s)$. Note that $\kappa_y(t_0)<0$.
	
	Since 
	\begin{equation*}
	y\in C_b(t_0)=\bigcup_{t\in\left(t_0,b\right)}B(\Gamma(t),r(t))\backslash \left(\overline{B(\Gamma(t_0),r(t_0))}\right),
	\end{equation*}
	there exist $\tilde{t}_1>t_0$ such that $y\in B\left(\Gamma(\tilde{t}_1,r(\tilde{t}_1))\right)$. This implies that $\kappa_y(\tilde{t}_1)>0$ and the continuity of $\kappa$ gives the existence of a point $t_1\in \R_+$ such that $y\in\overline{C_b(t_1)}$.
\end{proof}

\begin{proof}[Proof of Theorem~\ref{corollary of SMP}:]
	(See Figure~\ref{Fig1}) We argue by contradiction. Suppose that there exists a compact manifold $(M,g)$ and a non-constant harmonic map $\phi : M\longrightarrow N$  such that $\phi(M)\subset \mathcal{R}$. Let $\phi(x)=y\in\mathcal{R}$ be a point in $\phi(M)$. By definition of $\mathcal{R}$ and the fact that $\phi$ is non-constant, there exists $t_{0} \in [a,b]$ such that, 
	$$
	y \in \partial B(\Gamma(t_{0}),r(t_{0}))\cap \overline{C_b(t_0)},
	$$
	where $C_b(t)$ is defined by Equation~\eqref{connected component of gamma b}.
	Therefore, by Sampson's maximum principle, for every open neighborhood $U_{x}$ of $x$ in $M$, 
	$$
	\phi(U_{x})\nsubseteq B(\Gamma(t_{0}),r(t_{0}))
	$$
	That is, there exists $x_{1}\in U_{x}$ such that $\phi(x_{1})=y_{1}$ is mapped to the convex side of $B(\Gamma(t_0,r(t_0)))$. Moreover, as we assume that $\mathcal{R}\backslash B(\Gamma(t_{0}),r(t_{0}))$ is the union of two disjoint connected sets, $y_{1}\in C_b(t_0)$. Therefore, by Lemma~\ref{property of R}, there exists a $t_1>t_0$ such that 
	$$
	y_{1}\in\partial B(\Gamma(t_{1}),r(t_{1}))\cap \overline{C_b(t_1)}
	$$
	Again, by Sampson's maximum principle, every open neighborhood $U_{x_{1}}$ of $x_1$ in $M$ satisfies 
	$$
	\phi(U_{x_{1}})\nsubseteq B(\Gamma(t_{1}),r(t_{1})),
	$$ 
	that is, there exists $x_{2} \in U_{x_{1}}$ such that $\phi(x_{2})=:y_{2}\notin B(\Gamma(t_{1}),r(t_{1}))$, $y_{2}$ is in the convex side of $\partial B(\Gamma(t_{1}),r(t_{1}))$ and in the same connected component as $\Gamma(b)$. By induction, we construct sequences $\{t_{i}\}_{i\in \N}$ and $\{y_{i}\}_{i\in \N}$ such that for every $i\in \N$ we have $t_i>t_{i+1}$ and $y_i=\phi(x_i)\in \phi(M)\cap \overline{C_b(t_i)}$.
	
	Define $t^*=\sup\{t\in (a,b); \phi(M)\cap\overline{C_b(t)}\neq \emptyset\}$. We claim that $t^*=b$. Indeed, by definition of $t^*$ and the fact that $\{t_i\}$ is a bounded monotone sequence, we have $t_j\longrightarrow t^*$ and $y_j=\phi(x_j)\in \overline{C_b(t_j)}\subset\overline{C_b(t_{j-1})}\subset...\subset\overline{C_b(t_0)}$. Denoting by $y_{\infty}=\lim_{j\rightarrow\infty} y_j$, we have $y_\infty=\phi(x_\infty)\in\overline{C_b(t_0)}\cap...\cap\overline{C_b(t^*)}=\overline{C_b(t^*)}$ and if $t^*<b$ we can apply Sampson's maximum principle to get a time $t^{**}$ such that $t^*<t^{**}<b$ such that $\phi(M)\cap\overline{C_b(t^{**})}\neq \emptyset$ which contradicts the fact that $t^*$ is the supremum of those times. Hence $y_{\infty}\in\partial B(\Gamma(b),r(b))\cap\partial\mathcal{R}$ and using once more Sampson's maximum principle we obtain a point $\tilde{y}\in\phi(M)\setminus\mathcal{R}$.  
\end{proof}

	

\noindent This result  has very interesting applications.

\begin{figure}
	\begin{center}
	\includegraphics[width=8cm]{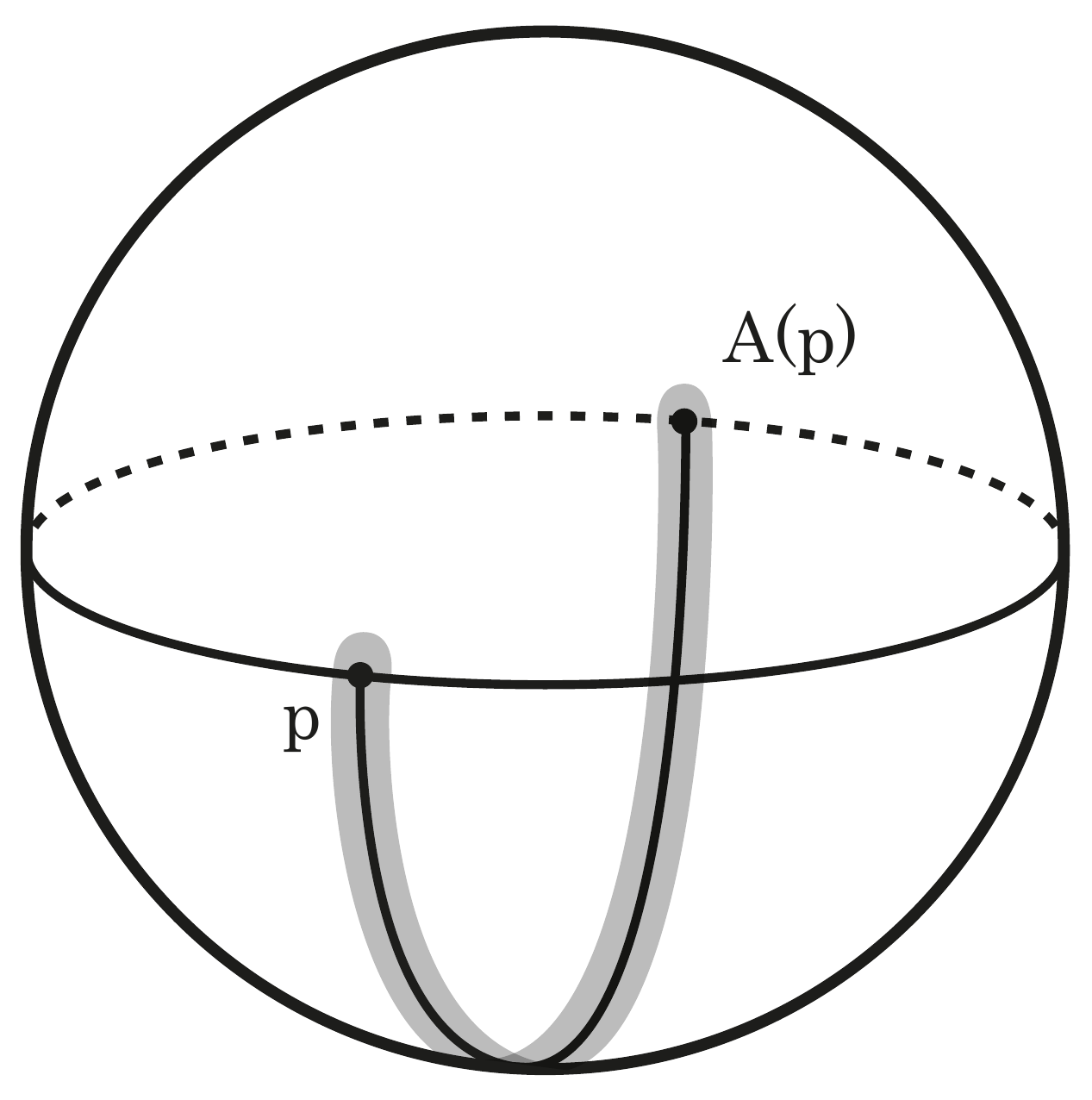}
	\caption{$S^2\backslash\left(\frac{1}{2}Eq\right)_{\epsilon}$}
	\label{Fig2}
	\end{center}
\end{figure}

\begin{eg}\label{S2 without half equator}
	Let $(N,h)=(S^2,\mathring{g})\subset (\R^3,g_{\text{euc}})$ and let $\Gamma: [0,1]\longrightarrow S^2$ be the unit speed geodesic such that $\Gamma(0)=(-1,0,0)$, $\Gamma(\frac{1}{2})=(0,0,1)$ and $\Gamma(1)=(1,0,0)$. For a given $\epsilon > 0$, let
	\begin{equation}
	\mathcal{R}:=\displaystyle\bigcup_{t\in[0,1]} B\left(\Gamma(t),\frac{\pi}{2}-\frac{\epsilon}{2}\right)
	\end{equation}
and note that $\mathcal{R}$ contains the complement of an $\epsilon$-neighborhood of the half equator containing the points $(0,0,-1)$, $(0,-1,0)$ and $(0,1,0)$, denoted by $\frac{1}{2}$Eq. We will write $\mathcal{R}$ as $S^2\backslash(\frac{1}{2}Eq)_{\epsilon}$. Theorem \ref{corollary of SMP}  implies that $S^2\backslash(\frac{1}{2}Eq)_{\epsilon}$ has property ($\star$).(See Figure~\ref{Fig2})
 
\end{eg} 
 
\begin{remark}
 	The above argument also works for spheres $S^n$ when we take out a closed half equator, denoted by $\frac{1}{2}S^{n-1}$. 
    This result appears in \cite{jost2012}, where is constructed a strictly convex function on each compact subset properly contained in  $S^n\backslash\frac{1}{2}S^{n-1}$.
\end{remark}

\begin{remark}
		
	The main idea of the  proof of Theorem \ref{corollary of SMP} is that $\partial \mathcal{R}$ is a barrier to the existence of non-constant harmonic maps. With the help of the SMP and the definition of $\mathcal{R}$ as a union of convex balls, we push the image of the harmonic map to this barrier. 
	
	One can make the theorem more general by simply changing $\partial \mathcal{R} \cap \partial B(\Gamma(b),r(b))$ to some more flexible barrier. For example, with a similar family of balls that we defined in Example~\ref{S2 without half equator}, one would prove that there are no non-constant harmonic maps whose image is contained in an open subset of $S^2\backslash\gamma$ where $\gamma$ is a connected curve connecting the antipodal points $(0,1,0), (0,-1,0)\in S^2$.
	
	Note that the classical method used to prove non-existence of such non-constant harmonic maps, i.e., the method of looking for a strictly convex function $f:\mathcal{R}\longrightarrow \R$ (in the geodesic sense), is not very flexible. Once one changes the boundary of $\mathcal{R}$ slightly, one can no longer guarantee that there exists a  strictly convex function $\tilde f:\mathcal{\tilde R} \longrightarrow \R$.

\end{remark}

\subsection{Barriers for the existence of harmonic maps}


Our proof of Theorem~\ref{corollary of SMP} gave the barrier for the existence of harmonic maps as the boundary of the union of a family of balls whose centers are given by a smooth embedded curve. 

By studying barriers to the existence of harmonic maps abstractly, we find that the result does not rely on the existence of such a family of balls along a curve. 

The region $\mathcal{R}$  defined in Example~\ref{S2 without half equator} contains no closed geodesics, and since geodesics are harmonic maps,  it is necessary  to avoid them for the non-existence of non-constant harmonic maps. Therefore the first condition we are looking for is one that avoids such curves.

\begin{condition} \label{i}
  For every $p\in\mathcal{R}$ and every $v\in T_{p}^{1}N$, the geodesic $\Gamma$ such that $\Gamma(0)=p$ and $\Gamma'(0)=v$ satisfies, for some $t_{v}\in [0,K]$, $\Gamma(t_{v})\in\partial\mathcal{R}$. Where $K$ is independent of $p$ and $v$.
\end{condition}

Condition~\ref{i} is a qualitative way of imposing the non-existence of closed geodesics in our region. But imposing the non-existence of closed geodesics in a given region is not enough to avoid existence of harmonic maps, therefore we need an extra condition to be stated later. But first, let us introduce some notation.

\begin{notation}
	For $p\in N$ and $v\in T_p N$, we denote by $\gamma_{p,v}:\R\longrightarrow N$ the geodesic that satisfies $\gamma_{p,v}(0)=p$ and $\gamma_{p,v}'(0)=v$.
\end{notation}

\noindent As an example, suppose we have a point in the image of some geodesic, namely $\Gamma_{p,v}(t_0)$. Take a tangent vector $\eta\in T_{\Gamma_{p,v}(t_0)} N$ such that $\eta\neq \Gamma_{p,v}'(t_0)$. If we consider another geodesic passing through ${\Gamma_{p,v}(t_0)}$ in the direction $\eta$ at $t=0$ , we write $\gamma_{\Gamma_{p,v}(t_0),\eta}(\cdot):\R\longrightarrow N$.

\begin{notation} 
   For $p\in N$, $\nu\in\partial B_{g}(p,r)$ with $r\ll 1$, let $t_{\max}(p,\nu)\in \R_{+}$ be the first time at which $\Gamma_{p,\nu}(t)\in\partial\mathcal{R}$.
\end{notation}

We impose the following condition on the region $\mathcal{R}$.

\begin{condition}\label{ii}
For every $p\in\mathcal{R}$ and every $\nu\in T_{p}^{1}N$, if $t<t_{\max}(p,\nu)$, then the set of directions $\left\{\eta\in T_{\Gamma_{p,\nu}(t)}^{1}N \,|\,  d(p,\gamma^{\mathcal{R}}_{\Gamma_{p,\nu}(t),\eta}(\epsilon)) < d(p,\Gamma_{p,\nu}(t))\right\}$ is properly contained in $T_{\Gamma_{p,\nu}(t)}^{1}$ and connected. Here $\gamma^{\mathcal{R}}$ is a geodesic in $\overline{\mathcal{R}}$ with $\epsilon\in\R_{+}$ sufficiently small (that is, a curve that minimizes the length in $\mathcal{R}$).
\end{condition}

In general, Condition~\ref{ii} is difficult to check unless $N$ is a symmetric space and one understands very well the geometry of the region $\mathcal{R}$. 

\begin{thm}\label{property star}
	Let (N,h) be a complete Riemannian manifold and suppose that for a region $\mathcal{R}\subset N$, $\mathcal{R}\setminus\partial\mathcal{R}$ is open and connected. In addition, suppose $\mathcal{R}$ satisfies the Conditions~\ref{i} and  \ref{ii} above. Then $\mathcal{R}$ has property ($\star$) of Definition~\ref{star}.
\end{thm}

\begin{proof}
	We argue by contradiction. Suppose that there exists a compact manifold $(M,g)$ and a non-constant harmonic map $\phi :M\longrightarrow N$  such that $\phi(M)\subset \mathcal{R}$. Therefore, there exist points $p=\phi(x),q \in\mathcal{R}$ such that $d(q,p) \ll 1$ and a geodesic $\Gamma$ with $\Gamma(0)=p$, $\Gamma^{\prime} (0)=v_{p,q}$, and $\Gamma(t_{p})=q$. By Condition~\ref{i} there exists $t_{v_{p,q}}\in \R_+$ such that $\Gamma$ intersects $\partial\mathcal{R}$. Define 
	$$
	Q_{0}:= \cup_{t\in\R_{+}}B(\Gamma(t),\delta), \,\text{with}\,\, \delta>\delta_{0}>0,
	$$
	and note that it satisfies the conditions of Theorem \ref{corollary of SMP}. 
	
	Therefore, when walking along $\Gamma$ with the ball $B(\Gamma(t),\delta)$, there exists $q_{1}\notin Q_{0}$ such that $q_{1}=\phi(y_{1})$, $y_{1}\in M$, and Condition~\ref{ii} tells us that 
	$$
	d_{\mathcal{R}}(q_{1},p)>d_{\mathcal{R}}(q,p).
	$$ 
	
	Let $\gamma^{\mathcal{R}}_{1}$ be a curve in $\overline{\mathcal{R}}$, connecting $p$ and $q_{1}$, of minimal length. Note that it may not be a geodesic in N itself. By Condition~\ref{ii}, there exists $\eta_{1}\in T^{1}_{q_{1}}N$ such that $\Gamma_{q_{1},\eta_{1}}$ increases the distance to $p$, that is, 
	$$
	d_{\mathcal{R}}(\Gamma_{q_{1},\eta_{1}}(t),p)>d_{\mathcal{R}}(q_{1},p)
	$$ 
	for sufficiently small t (as before, $\Gamma_{q_{1},\eta_{1}}$ hits the boundary at some point $t_{q_{1}}$ by Condition~\ref{i}). 
	
	Now let 
	$$Q_{1}:=\cup_{t\in\R_{+}}B(\Gamma_{q_{1},\eta_{1}}(t),\delta_{1}),\,\,\delta_{1}>\delta_{0}.$$ 
	
	By Theorem~\ref{corollary of SMP}, there exist $q_{2}\in Q_{1}$ such that $q_{2}=\phi(y_{2})$ and 
	$$
	d_{\mathcal{R}}(q_{2},p)>d_{\mathcal{R}}(q_{1},p)>d_{\mathcal{R}}(q,p).
	$$
	
	Inductively, we build a sequence of points $q_{1},q_{2},q_{3},...$ such that $q_{i}=\phi(y_{i})$ with $y_{i}\in M$ for all $i\in \N$ and 
	$$
	d_{\mathcal{R}}(q_{i},p)=d_{\mathcal{R}}(q_{i},\phi(x))>d_{\mathcal{R}}(q_{i-1},\phi(x)).
	$$
	
	Arguing like in the proof of \ref{corollary of SMP}, we obtain that $y_\infty=\phi(x_\infty)\in\phi(M)\cap\partial\mathcal{R}$ and therefore by Sampson's maximum principle there exists a point $\tilde{y}\in \phi(M)\setminus \mathcal{R}$.
\end{proof}

\begin{eg}
One can use the above theorem to find regions in symmetric spaces for which property $(\star)$ holds. If one has some information about the geometry of such a space, Conditions~\ref{i} and \ref{ii} may be easy to check. For instance, let $p\geq 1$ and $(\Sigma_p,g)$ be a compact surface of genus $p$ with a hyperbolic metric of constant curvature. Let $\gamma_1,...,\gamma_{2p}$ be smooth curves that generate $H_1(\Sigma_p,\R)$, the first homology group of $\Sigma_p$. For a given $\epsilon>0$, define 
\begin{equation}
\mathcal{R}:=\Sigma_p\backslash\left(\gamma_1\cup...\cup\gamma_{2p}\right)_{\epsilon}
\end{equation}
\noindent We see that $\mathcal{R}$ satisfies Conditions~\ref{i} and \ref{ii}, therefore it has property ($\star$). See Figure~\ref{Fig3}.

\begin{figure}
		\begin{picture}(100,200)
		\put(0,0){\includegraphics[width=1\linewidth]{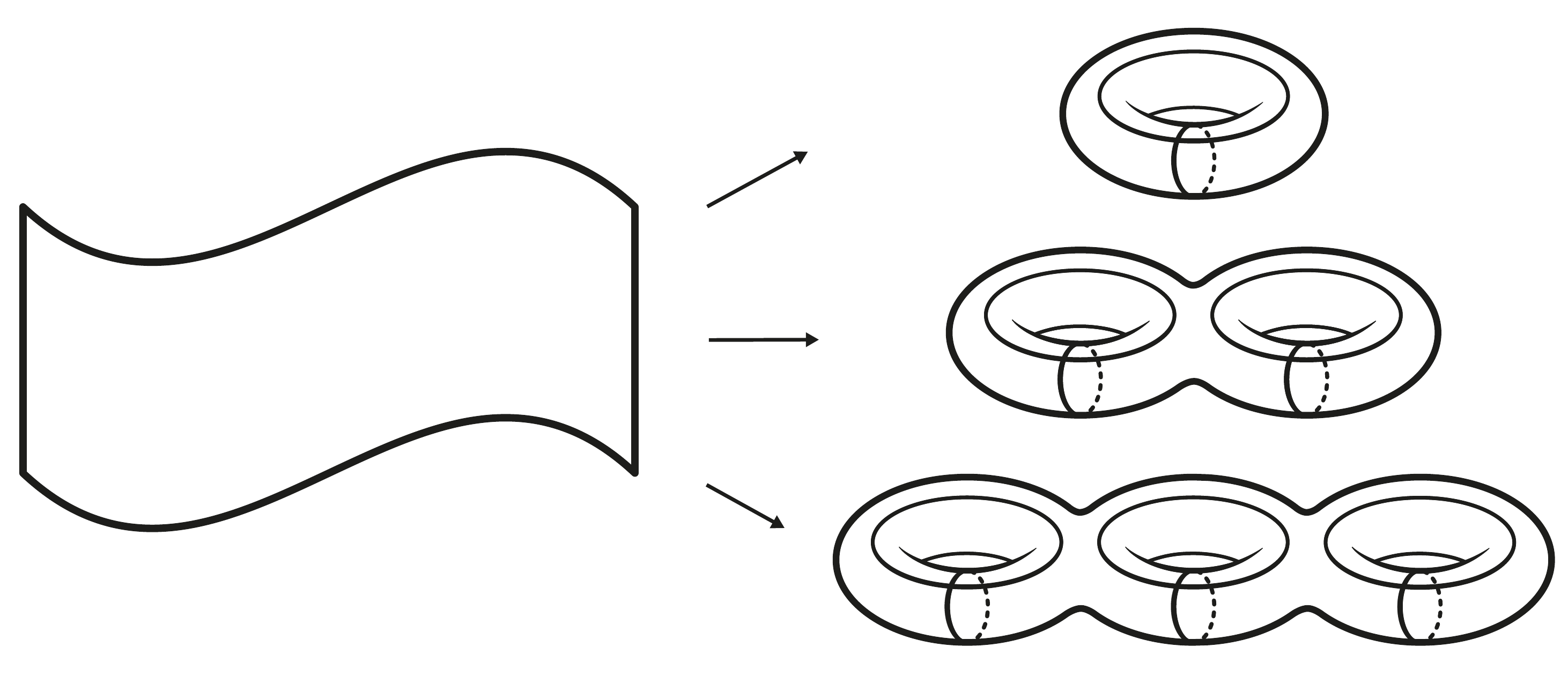}}
		\put(215,108){$\phi$}
		\put(195,83){harmonic}
		\put(125,45){$M$ compact}
		\put(125,25){$\partial M = \emptyset$}
		\end{picture}
	\caption{Example of barriers for the existence of harmonic maps in immersed surfaces in $\R^3$ of genus one, two and three.}
	\label{Fig3}
\end{figure}
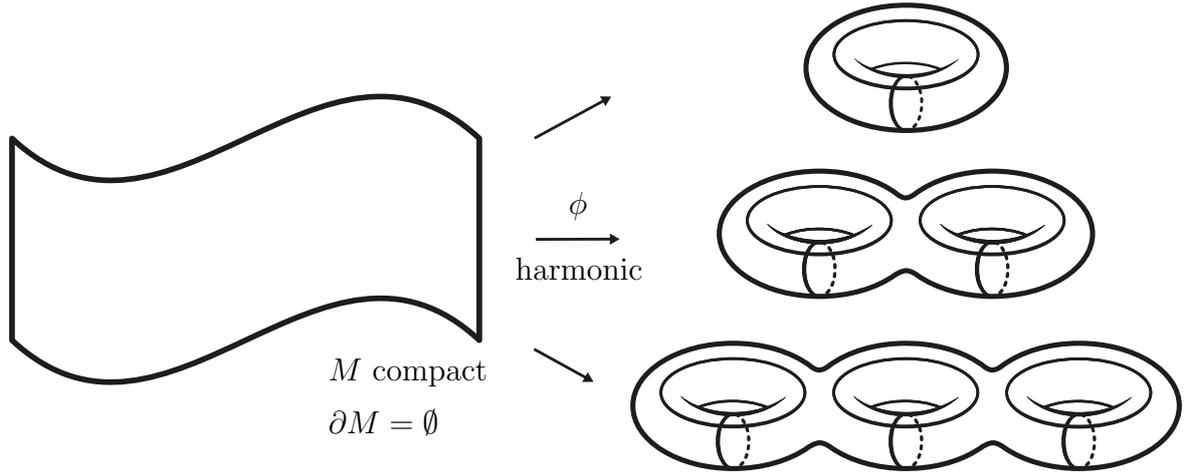

\noindent An alternative proof of this statement could be given simply by noting that a harmonic map defined on a compact manifold and taking values in a polygon with $4p$-sides in the Poincar\'e hyperbolic space such that its image is contained in the interior of this polygon must be constant. The example above is just to illustrate the relations between such property of harmonic maps and maximum principles.
\end{eg}

To use the previous results to obtain Bernstein-type theorems we have to overcome the compactness of $M$. We will be able to do so for the graph case and the harmonic map that we are taking into consideration is the Gauss map. We state that precisely in the following theorem.

\begin{thm}[Property ($\star$) for graphs]\label{property star for graphs}
	Let $M=\grafico(f)$ be a minimal submanifold, where $f:\R^p\longrightarrow \R^{n-p}$ is a smooth map with bounded slope, and let $\gamma:M\longrightarrow G^+_{p,n}$ be its normal Gauss map. Suppose $\gamma(M)\subset \mathcal{R}$, where $\mathcal{R}$ is a region in $G^+_{p,n}$ that has property ($\star$). Then the harmonic Gauss map is constant and $M$ is a plane.
	
\end{thm}

\begin{proof}	
	Let $C(M,\infty)$ be the tangent cone of $M$ at $\infty$ defined as follows. Considering the family  $f_t=\frac{1}{t}f(t.x)$ indexed on $R_+$ and the set $\left\{M_t=\grafico(f_t): t\in R_+\right\}$, it is known that $f_t$ satisfies the minimal surface equations and hence $\left\{M_t=\grafico(f_t): t\in R_+\right\}$ defines a family of minimal submanifolds in $\R^{n}$. Using the same argument as Fischer-Colbrie in \cite{fischer80}, we show that $\{f_t: t\in\R_+\}$ is an equicontinuous family on any compact subset of $\R^p$, thus Arzela-Ascoli implies that there exists a subsequence $f_{t_i}: i\in \N$ such that $\lim_{i\rightarrow +\infty}=+\infty$ and $\lim_{i\rightarrow +\infty}f_{t_i}=h$, where $h$ is a Lipschitz function and a weak solution to the minimal surface equation, and its graph $C(M,\infty):=\{(x,h(x)): x\in\R^n\}$ is a cone, which is called the tangent cone of $M$ at $\infty$.\\
	
	If the cone $C(M,\infty)$ is regular except in the origin, then the intersection of $C(M,\infty)$ and the unit sphere $S^{n-1}$ gives a $(p-1)$-dimensional embedded compact minimal submanifold of $S^{n-1}$, which is denoted by $M^{\prime}$ (\cite{simons68}). Since the normal Gauss map is a cone-like map, the image of the normal Gauss map of $M^{\prime}$ is still contained in the region $\mathcal{R}$. $M^{\prime}$ is compact and minimal in $S^{n-1}$, hence $M^{\prime}$ has parallel mean curvature and is compact in $\R^{n}$. Thus the fact that $\gamma(M^{\prime})\subset \mathcal{R}$ and the region $\mathcal{R}$ has property $(\star)$ implies that $M^{\prime}$ is a totally geodesic subsphere of $S^{n-1}$. Then Allard's regularity in \cite{allard72} implies that $f$ is an affine map and $M$ a $p$-plane.\\
	
	If the cone $C(M,\infty)$ has a singularity at a point $p_0\neq 0$, we consider a blow-up at $p_0$ and obtain a tangent cone of the tangent cone $C(M,\infty)$ at $p_0$, which is a product of the form $C\times\R$, and $C$ is a minimal cone of dimension $(p-1)$ in $\R^{n-1}$, which is also a graph and the slope is bounded everywhere. If $C$ is smooth except at the origin, we can proceed as above to prove that $C$ is an affine $(p-1)$-plane, which gives a contradiction. Otherwise, we iterate the above argument until we have obtained a minimal cone of dimension $1$, which cannot be singular once more giving a contradiction. Therefore $0\in R^{n}$ is the only possible singularity for $C(M,\infty)$ and hence $M$ is an affine $p$-plane.  
	
\end{proof}

\noindent In section 4, we will use this to study regions on the Grassmannian manifold where the harmonic Gauss map from a minimal submanifold may, or may not, be constant.

\section{The Geometry of Grassmannian Manifolds}
We first need to recall some facts about 
 Grassmannians as stated in Kozlov \cite{kozlov97} and Jost-Xin \cite{jost99}.

\subsection{Basic definitions}

Let $V^{n}$ be an $n$-dimensional vector space over $\R$ with inner product $\langle .,.\rangle$. One defines a $2^{n}$-dimensional algebra with respect to the exterior product $\wedge$ by
\begin{equation}
\Lambda(V)=\oplus_{i=0}^{\infty} \Lambda_{i}(V)
\end{equation}
such that $\Lambda_{0}(V)=\R$, $\Lambda_{1}(V)=V$ and $\Lambda_{i>n}(V)=0$. $\Lambda(V)$ is called the Grassmann algebra and the elements of $\Lambda_{p}(V)$ are called $p$-vectors. When $V=\R^{n}$, we denote $\Lambda_p(\R^n)$ simply by $\Lambda_p$.

\noindent Let $\{e_{i}\}_{i=1}^{n}$ be a basis for $V$ and $\lambda=\left(i_{1},...,i_{p}\right)\in \Lambda\left(n,p\right)=\{\left(i_{1},...i_{p}\right); 1\leq i_{1} < ... < i_{p}\leq n\}$. We denote by $e_{\lambda}$ the $p$-vector
\begin{equation}
e_{\lambda}:=e_{i_{1}}\wedge ... \wedge e_{i_{p}}
\end{equation}
$\{e_{\lambda}\}_{\lambda\in\Lambda\left(n,p\right)}$ is a basis for the $\binom{n}{p}$-dimensional vector space $\Lambda_{p}(V)\subset \Lambda(V)$. For any $w\in\Lambda_{p}(V)$,
\begin{equation}
w=\sum_{\lambda\in\Lambda(n,p)} w^{\lambda}e_{\lambda}
\end{equation}

\begin{defn}[Simple vector]
A $p$-vector $w\in\Lambda_{p}(V)$ is called simple if it can be represented as the exterior product of p elements of V, that is, $w=a f_{1}\wedge...\wedge f_{p}$, with $a\in\R$ and $f_{i}\in V$ for every $i\in\{1,...,p\}$.
\end{defn}

\begin{defn}[Scalar product of $p$-vectors]
	The scalar product of two $p$-vectors $w$ and $\tilde{w}$ is the number
	\begin{equation}
	\langle w,\tilde{w}\rangle:= \sum_{\lambda\in\Lambda(n,p)} w^{\lambda}\tilde{w}^\lambda
	\end{equation} 
\end{defn}

\noindent One can prove that this product has the properties of a scalar product and does not depend on the choice of the orthonormal basis. We also have that $\left|w\right|={\langle w,w\rangle}^{\frac{1}{2}}$.

\begin{defn}[Inner multiplication of $p$-vectors]
	An operation of inner multiplication is a bilinear map
	\begin{align}
	\llcorner:\Lambda_p(V)\times\Lambda_q(V) &\longrightarrow \Lambda_{p-q}(V), \hspace{1cm} p\geq q\geq 0\\
              (\omega,\xi) &\mapsto w\llcorner\xi
	\end{align}
	that for any $\omega\in\Lambda_p(V), \xi\in\Lambda_q(V)$, and $\varphi\in\Lambda_{p-q}(V)$ satisfies
	\begin{equation}
	\langle \omega\llcorner\xi,\varphi\rangle=\langle\omega,\xi\wedge\varphi\rangle
	\end{equation}
\end{defn}

\begin{defn}[The rank space]
	The rank space of a $p$-vector $w\in\Lambda_{p}(V)$ $\left(p\geq 1\right)$ is defined as
	\begin{equation}
	V_{w} = w \llcorner \Lambda_{p-1}(V)=\{e\in V; e=w\llcorner v, v\in\Lambda_{p-1}(V)\}\subset V
	\end{equation}
	That is, $V_{w}\subset V$ is the minimum subspace $W\subset V$ such that $w\in \Lambda_{p}(W)$.
\end{defn}

\begin{remark}
	A nonzero $w \in \Lambda_{p}(V)$ is simple if, and only if, $\dim V_{w}=p$.
\end{remark} 

\noindent Therefore there is a one-to-one correspondence between the set of oriented $p$-planes and the set of \emph{unit simple} $p$-vectors. This will be important when we define the so called Pl\"ucker embedding of the Grassmannian manifold.

\begin{defn}[Grasmannian manifold $G^{+}_{p,n}$]
	Let $\R^{n}$ be the n-dimensional Euclidean space. The set of all oriented $p$-subspaces (called $p$-planes) constitutes the Grassmannian manifold $G^{+}_{p,n}$, which is the irreducible symmetric space $SO(n)/\left(SO(p)\times SO(q)\right)$, where $q=n-p$.
\end{defn}

\noindent For each $p$-plane $V_{0}\in G^{+}_{p,n}$, consider  the open neighborhood $U(V_{0})$ of all oriented $p$-planes whose orthogonal projection into $V_{0}$ is one-to-one. That is, if $\{e_{i}\}_{i=1}^{p}$ is an orthonormal basis of $V_{0}$, and $\{n_{\alpha}\}_{\alpha=1}^{q}$ is an orthonormal basis of $V_{0}^{\perp}$, one can write, for $V\in U(V_{0})$,
\begin{equation}
e_{i}(V)=e_{i} + N_{i}
\end{equation}
where $N_{i}\in V_{0}^{\perp}$ and $e_{i}(V)\in V$ is such that $e_{i}=Pr_{V_{0}}(e_{i}(V))$.

\noindent Decomposing the vector $N_{i}$ using the basis $\{n_{\alpha}\}$ gives us a matrix $a=\left(a_{i}^{\alpha}\right)\in \R^{p.q}$ defined by
\begin{equation}
N_{i}=a_{i}^{\alpha}n_{\alpha}
\end{equation}
that can be regarded as local coordinates of the $p$-plane V in $U(V_{0})$.

\noindent Using the one-to-one correspondence between unit simple vectors and oriented $p$-planes $V\in G^{+}_{p,n}$ we define the Pl\"ucker embedding as
\begin{equation}
	\begin{split}
    \psi: G^{+}_{p,n}\longrightarrow \Lambda_{p}\\
    \psi(V):=\frac{\tilde{\psi}(V)}{\left|\psi(V)\right|},
    \end{split}
\end{equation}
where
\begin{equation}
\tilde{\psi}(V):=(e_{1}+N_{1})\wedge...\wedge(e_{p}+N_{p}),\hspace{0.2cm}.
\end{equation}

Via the Pl\"ucker embedding, $G^{+}_{p,n}$ can be viewed as a submanifold of the Euclidean space $\R^{\binom{n}{p}}$. The restriction of the Euclidean inner product denoted by $\omega:G^{+}_{p,n}\times G^{+}_{p,n}\longrightarrow \R$, is
\begin{equation}
\omega(P,Q)=\frac{\langle\psi(P),\psi(Q)\rangle}{\langle\psi(P),\psi(P)\rangle^{\frac{1}{2}}\langle\psi(Q),\psi(Q)\rangle^{\frac{1}{2}}}.
\end{equation}
If $\{e_1,..., e_p\}$ is an oriented orthonormal basis of P and $\{f_1,..., f_p\}$ is an oriented orthonormal basis of Q, then
\begin{equation}
w(P,Q)=\langle e_1\wedge...\wedge e_p,f_1\wedge...\wedge f_p\rangle = \det W,
\end{equation}
where $W=(\langle e_i, f_j\rangle)$.

\begin{remark}
	Note that $\psi(G^{+}_{p,n})=K_{p}\cap S^{\binom{n}{p} - 1}\subset \R^{\binom{n}{p}}\cong\Lambda_{p}$.
\end{remark}
\begin{remark}
	Let $w\in \psi(G^{+}_{p,n})$ and $\{e_{i}\}_{i=1}^{p}$, $\{n_{\alpha}\}_{\alpha=1}^{q}$ be an orthonormal basis for $V_{w}$ and $V_{w}^{\perp}$, respectively. The system of vectors $\{\eta_{i\alpha}\}_{i=1,\alpha=1}^{p\hspace{0.2cm},\hspace{0.2cm}q}$ where
	\begin{equation}
	\eta_{i\alpha} = e_{1}\wedge...\wedge e_{i-1}\wedge n_{\alpha} \wedge e_{i+1} \wedge...\wedge e_{p},
	\end{equation}
	for $i\in\{1,...,p\}$ and $\alpha\in\{1,...,q\}$ is an orthonormal basis of the tangent space $T_{w}(\psi(G^{+}_{p,n}))$.
\end{remark}

We will write $G^+_{p,n}$ instead of $\psi(G^{+}_{p,n})$ for the image of the Grassmannian manifold under the Pl\"ucker embedding.

\noindent Let $(w,X)\in TG^{+}_{p,n}$ be an element of the tangent bundle and $\{\eta_{i\alpha}\}_{i=1,...,p}^{\alpha=1,...,q}$ a basis for $T_{w}G^{+}_{p,n}$ like above. There exist $m_{i}\in V_{w}^{\perp}$ such that
\begin{equation}
X = m_{1}\wedge e_{2}\wedge...\wedge e_{p}+...+e_{1}\wedge...\wedge e_{p-1}\wedge m_{p}.
\end{equation}
These $m_{i}$ are not necessarily pairwise orthogonal.

\begin{thm}[Kozlov]
	Let $w\in G^{+}_{p,n}$ and $X\in T_{w}G^{+}_{p,n}$, $X\neq 0$. Then there exists an orthonormal basis $\{e_{i}\}_{i=1}^{p}$ in $V_{w}$ and a system $\{m_{i}\}_{i=1}^{r}$, with $1\leq r\leq min\{p,q\}$, of non-zero pairwise orthogonal vectors in $V_{w}^{\perp}$, such that
\begin{align}
w =& e_{1}\wedge...\wedge e_{p}, \label{local expression of w}\\
X =& \left(m_{1}\wedge e_{2}\wedge...\wedge e_{r}+...+e_{1}\wedge...\wedge e_{r-1}\wedge m_{r}\right)\wedge\left(e_{r+1}\wedge...\wedge e_{p}\right).\label{local expression of X}
\end{align}
\end{thm}

\begin{notation}
	Based on Equation~\eqref{local expression of w} and Equation~\eqref{local expression of X} above, we write 
\begin{align}
&\lambda^{i}:= \left|m_{i}\right|, \hspace{0.2cm} n_{i}:= \frac{m_{i}}{\lambda^{i}},\\
&e_{i}(s):= e_{i}\cos(s)+n_{i}\sin(s), \hspace{0.2cm} n_{i}(s):= -e_{i}\sin(s)+n_{i}\cos(s), \\
&X_{0}:= e_{r+1}\wedge...\wedge e_{p}.
\end{align}
Note that $e_{i}^{\prime}(s)=n_{i}(s)$ and $n_{i}^{\prime}= -e_{i}(s)$, 
\end{notation}

\subsection{Closed geodesics in $G^{+}_{p,n}$.}

Consider the curve
\begin{equation}
w(t):=e_{1}(\lambda^{1}t)\wedge...\wedge e_{r}(\lambda^{r}t)\wedge X_{0}.
\end{equation}

\begin{thm}
	The curve $w(t)$ is a normal geodesic in the manifold $G^{+}_{p,n}$ and
	\begin{equation}
	\exp_{w}X = w(1)
	\end{equation}
\end{thm}
\begin{proof}[Proof(Sketch)]
	For each $t\in\R$, $\{e_{i}(\lambda^{i}t),n_{i}(\lambda^{i}t)\}_{i=1}^{r}$ are pairwise orthogonal vectors and
	
	\begin{align}
	e_{i}^{\prime}(\lambda^{i}t)=& \lambda^{i}n_{i}(\lambda^{i}t), \hspace{0.2cm} n_{i}^{\prime}(\lambda^{i}t)= -\lambda^{i}e_{i}(\lambda^{i}t), \\
	X(t):=& w^{\prime}(t) = \sum_{i=1}^{r}\lambda^{i}n_{i}(\lambda^{i}t)\wedge (w(t)\llcorner e_{i}(\lambda^{i}t)), \label{w prime} \\ 
	w(0)=w&, \hspace{0.2cm} w^{\prime}(0)=X, \hspace{0.2cm} \left|X(t)\right|=\left[\sum_{i=1}^{r}(\lambda^{i})^{2}\right]^{\frac{1}{2}}.
	\end{align}
Therefore $w(t)$ is parametrized proportionally to the arc-length. Moreover, taking an additional derivative, we have
\begin{equation}
w^{\prime\prime}(t) = -\left|X\right|^{2}w(t) + \xi(t),
\end{equation} 
where $\xi(t)$ is a term like Equation~\eqref{w prime} with  $n_{i}$ replacing the $e_{i}$ and therefore orthogonal to the tangent plane $T_{w(t)}G^{+}_{p,n}$.
\end{proof}

\noindent We will denote the above geodesic by
\begin{equation}\label{geodesic eq}
w_{X}(t):= e_{1}(\lambda^{1}t)\wedge...\wedge e_{r}(\lambda^{r}t)\wedge X_{0},
\end{equation}

\noindent where $X=(\lambda^{1}n_{1}\wedge e_{2}\wedge...\wedge e_{r} + \lambda^{2}e_{1}\wedge n_{2}\wedge...\wedge e_{r}+...+\lambda^{r}e_{1}\wedge...\wedge e_{r-1}\wedge e_{r})\wedge X_{0}$.

\noindent Thus a geodesic in $G^{+}_{p,n}$ between two $p$-planes is simply obtained by rotating one into the other in the Euclidean space, by rotating corresponding basis vectors. As an example the 2-plane spanned by $e_{1}$, $e_{2}$ in $\R^{4}$: one tangent direction in $G^{+}_{2,4}$ would be to move $e_{1}$ into $e_{3}$ and keep $e_{2}$ fixed, and the other tangent direction  would be to move $e_{2}$ into $e_{4}$ and keep $e_{1}$ fixed. In other words, we are taking $X,Y\in T_{w}G^{+}_{2,4}$, for some $w\in G^{+}_{2,4}$, where $X=e_{3}\wedge e_{2}$ and $Y=e_{1}\wedge e_{4}$ and the respective geodesics are given by equation $\eqref{geodesic eq}$. We can also consider $Z\in T_{w}G^{+}_{2,4}$ where $Z := (\frac{1}{\sqrt{2}}e_{3}\wedge e_{2} + \frac{1}{\sqrt{2}}e_{1}\wedge e_{4})$ and the respective geodesic $w_{Z}(t)$ given by equation $\eqref{geodesic eq}$ is obtained when we simultaneously rotate $e_{1}$ into $e_{3}$ and $e_{2}$ into $e_{4}$. This is a geometric picture that helps to understand our subsequent constructions. Later we compute the length of these different types of geodesics in a general Grassmannian.

\begin{thm}
	For $G^{+}_{p,n}$, we have
	\begin{equation}
	\mathrm{diam}(G^{+}_{p,n}) = \max\big\{\pi,\sqrt{r_{0}}\pi/2\big\},
	\end{equation}
where $r_{0}=min\{p,q\}$. More explicitly,
\begin{align}
r_{0} < 4 \Rightarrow& \hspace{0.2cm}\mathrm{diam}(G^{+}_{p,n}) = \pi, \\
r_{0} \geq 4 \Rightarrow& \hspace{0.2cm}\mathrm{diam}(G^{+}_{p,n}) = \sqrt{r_{0}}\,\,\frac{\pi}{2}. 
\end{align}
\end{thm}

For the next examples, we have to understand  the closed geodesics in $G^{+}_{p,n}$. 

\begin{thm}\label{closed geodesic in grassmannians}
	For geodesics in $G^{+}_{p,n}$ with parametrization \eqref{geodesic eq}, we have that
    \begin{equation}
    w_{X}(t_{1})=w_{X}(t_{2}) \iff \exists k\in\Z^{r}; \hspace{0.2cm}\lambda^{i}(t_{2}-t_{1})=\pi k^{i}, \hspace{0.1cm}\text{and}\hspace{0.1cm} \sum_{i=1}^{r}k^{i}\equiv 0 \mod 2. 
    \end{equation}
\end{thm}

\begin{rmk}
	    In particular, every geodesic loop in $G^{+}_{p,n}$ is a closed geodesic.
\end{rmk}

\noindent This assertion can be used to compute lengths of closed geodesics in $G^{+}_{p,n}$, because
\begin{equation}
L(w_X) = \left|X\right| = \left[\sum_{i=1}^{r}(\lambda^{i})^{2}\right]^{\frac{1}{2}}.
\end{equation}

\subsection{Geodesically convex sets in $G^{+}_{p,n}$.}

As we have mentioned in the previous section, we are interested in the convex sets in a Grassmannian manifold.

Let $w\in G^{+}_{p,n}$ and $X\in T_{w}G^{+}_{p,n}$ a unit tangent vector. We know that there exists an orthonormal basis $\{e_{i},n_\alpha\}_{i=1,...,p}^{\alpha=1,...,q}$ of $\R^n$ and a number $r\leq \min(p,q)$ such that $w=\test \{e_i\}$, $n=p+q$  and 
\begin{equation}\label{equation for X}
X=\left(\lambda^{1}n_{1}\wedge e_{2}\wedge...\wedge e_{r} +...+ \lambda^{r}e_{1}\wedge ... \wedge e_{r-1}\wedge n_{r}\right)\wedge X_{0},
\end{equation} 
and $\left|X\right|=\left(\sum_{\alpha=1}^{r}\left|\lambda^{\alpha}\right|^{2}\right)^{\frac{1}{2}}=1$.

\begin{defn}
	Let $X\in T_{w}G^{+}_{p,n}$ be given by Equation~\eqref{equation for X}. Define
	\begin{equation}
	t_{X}:=\frac{\pi}{2\left(\left|\lambda_{\alpha^{\prime}}\right| + \left|\lambda_{\beta^{\prime}}\right|\right)},
	\end{equation}
	where $\left|\lambda_{\alpha^{\prime}}\right|:=max\{\lambda^{\alpha}\}$ and $\left|\lambda_{\beta^{\prime}}\right|:=max\{\lambda^{\alpha}$; 
	$\lambda^{\alpha}\neq\lambda_{\alpha^{\prime}}\}$.
\end{defn}

\begin{defn}
    Let $w\in G^{+}_{p,n}$ and $X\in T_{w}G^{+}_{p,n}$ a unit tangent vector like in Equation~\eqref{equation for X}. Define
	$$B_{G}(w):=\{w_X(t)\in G^{+}_{p,n}; 0\leq t\leq t_{X}\}.
	$$
\end{defn}

\begin{thm}[Jost-Xin]\label{Jost-Xin}
	The set $B_{G}(w)$ is a (geodesically) convex set and contains the largest geodesic ball centered at $w$.
\end{thm}

\noindent Now we are going to present some explicit examples of geodesics in $G^{+}_{p,n}$ and compute at which time they  intersect the region we have defined above. Those examples clarify how geodesics behave in Grassmannian manifolds,  and  thereby the geometry of $B_{G}(w)$.

\begin{eg}\label{X1 type vector}
 Let $w=e_{1}\wedge...\wedge e_{p}\in G^{+}_{p,n}$ and denote $X_{1}=n_{1}\wedge e_{2}\wedge...\wedge e_{p}$. From the notation above $\lambda^{1}=1$ and $\lambda^{i}=0$ for every $i\neq 1$. Moreover $\left|X_{1}\right|=1$ and $t_{X_{1}}=\frac{\pi}{2}$. Thus the length of this closed geodesic is $2\pi$. Moreover, denoting by $\left|L(w_{X}(t))\right|$ the length of the segment connecting $w_{X}(0)$ and $w_{X}(t)$, we have $\left|L(w_{X_{1}}(t_{X_{1}}))\right|=\frac{\pi}{2}$. Clearly, $w_{X_{1}}(\pi)=-w$ and, as $G^{+}_{p,n}$ is a symmetric space, $w_{X_{1}}(-\epsilon)=w_{-X_{1}}(\epsilon)$, so $w_{X_{1}}:\R\longrightarrow G^{+}_{p,n}$ is a geodesic such that $\left. \gamma \right|_{\left[-\frac{\pi}{2},\frac{\pi}{2}\right]}\subset \overline{B_{G}(w)}$.
  
\noindent It is known that $G^{+}_{p,n}$ is a submanifold of $S^{\binom{n}{p}-1}\subset \R^{\binom{n}{p}}$ with the induced metric of the sphere. So in the direction $X_{1}$, $\overline{B_{G}(w)}$ contains half of a great circle connecting two antipodal points.
\end{eg}

\begin{eg}\label{X2 type vector}
	Consider $w=e_{1}\wedge...\wedge e_{p}\in G^{+}_{p,n}$ and $X_{2}=\frac{1}{\sqrt{2}}(n_{1}\wedge e_{2} + e_{1}\wedge n_{2})\wedge e_{3}\wedge...\wedge e_{p}$. Therefore $X_{2}\in T_{w}G^{+}_{p,n}$ with $\left|X_{2}\right|=1$. Moreover in the definition of $B_{G}(w)$:
	\begin{equation}
	t_{X_{2}}=\frac{\pi}{2(\frac{1}{\sqrt{2}}+\frac{1}{\sqrt{2}})}=\frac{\sqrt{2}\pi}{4}=\frac{\sqrt{2}}{2}\frac{\pi}{2}<\frac{\pi}{2}.
	\end{equation}
Thus, to make $t_{X}$ as small as possible, we need an equidistribution of the eigenvalues $\lambda^{\alpha}$ of $X$ satisfying $\sum_{\alpha}\left|\lambda^{\alpha}\right|=1$ and we maximize $\lambda_{\max}+\lambda_{\max^{\prime}}$. In this case,  $w_{X_{2}}\left(\frac{\sqrt{2}\pi}{4}\right)\in \partial B_{G}(w)$ is the  point in the boundary closest to $w$. Also, by Theorem~\ref{closed geodesic in grassmannians}
	\begin{equation}
	w_{X_2}(t) = w_{X_2}(0) = w \iff \exists k^1,k^2\in\Z^2;\hspace{0.2cm} \lambda^i(t-0)=\pi.k^i\hspace{0.2cm}\text{and}\hspace{0.2cm}  k^1 + k^2 = 0 \hspace{0.2cm}(\text{mod}\hspace{0.2cm}2)
	\end{equation}
	but
	\begin{equation}
	\lambda^1=\lambda^2=1/\sqrt{2}\hspace{0.1cm},\hspace{0.2cm}\text{therefore}\hspace{0.2cm} t=\sqrt{2}k^1\pi=\sqrt{2}k^2\pi.
	\end{equation}
	Then $k^1=k^2=1$ solves the equation and tells us that $w_{X_2}(\sqrt{2}\pi)=w$, and, as $w_{X_2}$ is parametrized by arc length, 
	\begin{equation}
	L(w_{X_{2}})=\sqrt{2}\pi=4t_{X_2}.
	\end{equation}
	By \cite{kozlov97}, this is the smallest non-trivial closed geodesic in $G^{+}_{p,n}$, $r_{0}>1$.

    \noindent It  will be crucial below that the equation for the geodesic $w_{X_{2}}$ is
    \begin{equation}\label{wX2}
    w_{X_{2}}(t)=\left(e_1 \cos\left(\frac{t}{\sqrt{2}}\right)+ n_1 \sin\left(\frac{t}{\sqrt{2}}\right)\right)\wedge\left(e_2 \cos\left(\frac{t}{\sqrt{2}}\right)+ n_2 \sin\left(\frac{t}{\sqrt{2}}\right)\right)\wedge \left( e_3\wedge...\wedge e_p\right),
    \end{equation}
    so we compute, for any $t\in\R$,
    
    \begin{align}\label{w inner wX2}
    \begin{split}
    \langle w_{X_2}(t),w\rangle = &\,\big\langle \big(e_1 \cos(t/\sqrt{2})+ n_1 \sin(t/\sqrt{2})\big)\\
    &\wedge\big(e_2 \cos(t/\sqrt{2})+ n_2 \sin(t/\sqrt{2})\big),(e_1\wedge e_2)\big\rangle\cdot\\
    &\big\langle(e_3\wedge...\wedge e_p),(e_3\wedge...\wedge e_p)\big\rangle\\
    = &\,\big\langle e_1 \cos\left(t/\sqrt{2}\right),e_1\big\rangle\cdot \big\langle e_2 \cos\left(t/\sqrt{2}\right),e_2\big\rangle = \cos^2\left(\frac{t}{\sqrt{2}}\right)\geq 0
    \end{split}
   \end{align}
    
    Therefore, using the fact that $G^{+}_{p,n}\subset S^{\binom{n}{p}-1}$ and denoting the antipodal point in $S^{\binom{n}{p}}$ by $-w$, the geodesic $w_{X_2}$ never leaves the hemisphere of $S^{\binom{n}{p}}$ that has $w$ as a pole. In particular $w_{X_2}(t\in\R_+)$ does not intersect $B_G(-w)$.
\end{eg}

\begin{eg}\label{Xr0 type vector}
	Consider $w=e_{1}\wedge...\wedge e_{p}\in G^{+}_{p,n}$ and $X_{r_{0}}=\frac{1}{\sqrt{r_{0}}}(n_{1}\wedge e_{2}\wedge...\wedge e_{r_{0}} +...+ e_{1}\wedge...\wedge e_{r_{0}-1}\wedge n_{r_{r_{0}}})\wedge e_{r_{0}+1}\wedge...\wedge e_{p}$, so that $X_{r_{0}}\in T_{w}G^{+}_{p,n}$ with $\left|X_{r_{0}}\right|=1$. Therefore
	\begin{equation}
	t_{X_{r_{0}}}=\frac{\pi}{2(\frac{1}{\sqrt{r_{0}}}+\frac{1}{\sqrt{r_{0}}})}=\frac{\sqrt{r_{0}}\pi}{4},
	\end{equation}
	and the geodesic $w_{X_{r_{0}}}$ lies inside $B_{G}(w)$. As before, $G^{+}_{p,n}$ is a symmetric space, thus the curve $w_{X_{r_{0}}}:(-\frac{\sqrt{r_{0}}\pi}{4},\frac{\sqrt{r_{0}}\pi}{4})\longrightarrow G^{+}_{p,n}$ is a geodesic with $w_{-X_{r_{0}}}(-\frac{\sqrt{r_{0}}\pi}{4}) = w_{X_{r_{0}}}(\frac{\sqrt{r_{0}}\pi}{4})$ and 
	\begin{align}\label{closed geodesic condition}
	\begin{split}
	w_{X_{r_0}}(t) = w_{X_{r_0}}(0) = w \iff& \exists k^1,...,k^{r_0}\in\Z^{r_0};\hspace{0.2cm} \lambda^i(t-0)=\pi.k^i\\
	&\text{and}\hspace{0.2cm}  \sum_{r_0 = 1}^{r_0}k^i = 0 \hspace{0.2cm}(\text{mod}\hspace{0.2cm}2)
	\end{split}
	\end{align}
	thus $\lambda^i=1/\sqrt{r_0}\hspace{0.2cm}$ and $t=\sqrt{r_0}k^i\pi$ for every $i\in\{1,...r_0\}$.
	
	\noindent Consequently, if $r_{0}$ is an even number, then $k^i=1$ solve the equation and tell us that $w_{X_{r_0}}(\sqrt{r_0}\pi)=w$, and, as $w_{X_{r_0}}$ is parametrized by arc length, $\left|L(w_{X_{r_{0}}})\right|=\sqrt{r_{0}}\pi\hspace{0.1cm}(=4t_{r_0})$; if $r_0$ is an odd number, then $k^i=2$ solve the equation and tell us that $w_{X_{r_0}}(\sqrt{r_0}\pi)=w$, and, as $w_{X_{r_0}}$ is parametrized by arc length, $\left|L(w_{X_{r_{0}}})\right|=2\sqrt{r_{0}}\pi\hspace{0.1cm}(=8t_{r_0})$.
	
\noindent	The equation for the geodesic $w_{X_{r_0}}$ is
	
	\begin{align}\label{wXr0}
	\begin{split}
	w_{X_{r_0}}(t)=&\left(e_1 \cos\left(\frac{t}{\sqrt{r_0}}\right)+ n_1 \sin\left(\frac{t}{\sqrt{r_0}}\right)\right)\wedge\left(e_2 \cos\left(\frac{t}{\sqrt{r_0}}\right)+ n_2 \sin\left(\frac{t}{\sqrt{r_0}}\right)\right)\wedge\\
	& ...\wedge \left(e_{r_0} \cos\left(\frac{t}{\sqrt{r_0}}\right)+ n_{r_0} \sin\left(\frac{t}{\sqrt{r_0}}\right)\right)\wedge \left(e_{r_0 + 1}\wedge...\wedge e_p\right)
	\end{split}
	\end{align}
	so we compute for any $t\in\R_+$
	\begin{align}\label{w inner wXr0}
	\begin{split}
	\langle w_{X_{r_0}}(t),w\rangle = &\quad\langle \big(e_1 \cos(t/\sqrt{r_0})+ n_1 \sin(t/\sqrt{r_0})\big),e_1\rangle\cdot\\
	& \quad\langle \big(e_2 \cos(t/\sqrt{r_0})+ n_2 \sin(t/\sqrt{r_0})\big),e_2\rangle\cdots\\
	& \quad\langle \big(e_i \cos(t/\sqrt{r_0})+ n_i \sin(t/\sqrt{r_0})\big),e_i\rangle \cdots\\
	& \quad\langle \big(e_{r_0} \cos(t/\sqrt{r_0})+ n_{r_0} \sin(t/\sqrt{r_0})\big),\rangle\cdot\\
	&\quad\langle(e_{r_0 +1}\wedge...\wedge e_p),(e_{r_0+1}\wedge...\wedge e_p)\rangle\\
	= &\quad\langle e_1 \cos(t/\sqrt{r_0}),e_1\rangle\cdots \langle e_{r_0} \cos(t/\sqrt{r_0}),e_{r_0}\rangle = \cos^{r_0}(t/\sqrt{r_0})
	\end{split}
	\end{align}
	
\noindent Note that if $r_0$ is  odd, then $w_{X_{r_{0}}}(0)$ and $w_{X_{r_{0}}}(\sqrt{r_{0}}\pi)$ are antipodal points in $S^{\binom{n}{p}-1}$. Also, $\langle w_{X_{r_0}}(t),w\rangle \geq 0$ and $w_{X_0}$ does not intersect $B_G(-w)$ if and only if $r_0$ is even.
\end{eg}

\begin{notation}[Eigendirections on the tangent space]
The elements of the tangent space $Y\in T_{w}G^{+}_{p,n}$ given by the rotation of $1\leq k\leq r_0$ basis vectors into $k$ normal vectors with the same speed (i.e., an integer multiple of $\lambda_k$) will be called vectors of type k. For example, in Example~\ref{X1 type vector}, $X_1\in T_{w}G^+_{p,n}$ is a vector of type 1. In Example~\ref{X2 type vector}, $X_2\in T_{w}G^{+}_{p,n}$ is a vector of type 2.
\end{notation}

\section{Main Results}

In this section we construct regions in Grassmannian manifolds that have property $(\star)$ for graphs and use them to prove Bernstein-type theorems. First, we need to explain the relation of  the geometry of Gauss maps with the slope of  graphs. To show the gist of our method, we first explain how Moser's Bernstein theorem for codimension 1 follows from our method. Then we proceed to derive an analogous result in codimension 2. 

\subsection{The slope of a graph}

In our main theorem the concept of slope of a graph, defined by 
\begin{equation}
\Delta_{f}(x):=\big\{\det\big(\delta_{\alpha\beta} + \sum_{i}f^{i}_{x_{\alpha}}(x)f^{i}_{x_{\beta}}(x)\big)\big\}^{\frac{1}{2}},
\end{equation}
and its relation with harmonic Gauss maps plays an important role which we now describe.

Let $f^{\alpha}(x^1,...,x^n):\R^n\longrightarrow \R$ be smooth functions with $\alpha=1,...,m$ and let 
$$
M=\big(x^1,...,x^n,f^1(x^1,...,x^n),...,f^m(x^1,...,x^n)\big)\subset \R^{n+m},\,\ie,\, M=\grafico(f)\subset \R^{n+m}.
$$ 
Denote by $F:\R^n\longrightarrow M$ the diffeomorphism $x\mapsto (x,f(x))$, where $f=(f^1,...,f^n)$. $M$ is an $n$-dimensional submanifold of the Euclidean space $\R^{n+m}$ with the induced metric $g=g_{ij}dx^{i}dx^{j}$. With a basis $\{e_i,\eta_{\alpha}\}$ for $\R^{n+m}$, where $\{e_i\}$ spans $\R^n$, we can write this metric as
\begin{equation}
g_{ij} = \langle F_{*}\frac{\partial}{\partial x^i},F_{*}\frac{\partial}{\partial x^j}\rangle=\langle e_i + \frac{\partial f^{\alpha}}{\partial x^i}\eta_{\alpha}, e_j + \frac{\partial f^{\alpha}}{\partial x^j}\eta_{\alpha}\rangle = \delta_{ij} + \frac{\partial f^{\alpha}}{\partial x^i}\frac{\partial f^{\alpha}}{\partial x^j},
\end{equation} 
\noindent so denoting by $Df:=(\frac{\partial f^{\alpha}}{\partial x^i})$, a $(n\times m)$-matrix, \begin{equation}(g_{ij})=I_{n} + Df(Df)^{T}.\end{equation}

But also, taking $w,v\in G_{n,n+m}$ such that $w=e_1\wedge...\wedge e_n$, $v=e_{n+1}\wedge e_2\wedge...\wedge e_n$ and using the fact that
\begin{equation}
\psi\circ\gamma=(e_1 + \frac{\partial f^{\alpha}}{\partial x^1}\eta_{\alpha})\wedge...\wedge(e_n + \frac{\partial f^{\alpha}}{\partial x^n}\eta_{\alpha}).
\end{equation}
where $\psi:G^{+}_{n,n+m}\longrightarrow \R^{\binom{n}{n+m}}$ is the Pl\"ucker embedding,  we see that
\begin{equation}\label{the slope of the gauss map}
\omega(\gamma,w) = \frac{\langle\psi\circ\gamma,w\rangle}{\langle\psi\circ\gamma,\psi\circ\gamma\rangle^{\frac{1}{2}}\langle w,w\rangle^{\frac{1}{2}}} =  \Delta_{f}^{-1},
\end{equation}
\begin{equation}
\omega(\gamma,v) = \frac{\langle\psi\circ\gamma,v\rangle}{\langle\psi\circ\gamma,\psi\circ\gamma\rangle^{\frac{1}{2}}\langle v,v\rangle^{\frac{1}{2}}} =  \frac{\partial f^1}{\partial x^1}\Delta_{f}^{-1}.
\end{equation}
Where $\omega(.,.)$ is the restriction of the Euclidean inner product to the image of the Pl\"ucker embedding, as defined in section 3.

The above equations tell us that Moser's condition of a bounded slope of $\grafico(f)$, $\Delta_{f}<\infty$, in terms of the geometry of Grassmannians says that there exists a point $w\in \gamma(M)$ such that $\gamma(M)\subset G^{+}_{n,n+m}\cap \Hem_{w}^{\binom{n}{n+m}-1}$, where $\Hem_{w}^{\binom{n}{n+m}-1}$ is the hemisphere in $S^{\binom{n}{n+m}-1}$ centered at $w$. Obviously, as in the codimension 1 case, $G^{+}_{n,n+1}=S^{n}$ and the hemisphere of a point is a convex set, we have a proof for Moser's Bernstein-type theorem in \cite{moser61}.

\subsection{Bernstein-type theorems for codimension 2}

Here we want to use the geometry of Grassmannians to construct regions in $G^{+}_{p,p+2}$ satisfying Conditions~\ref{i} and \ref{ii} in Theorem \ref{property star}, to conclude  that these regions have property ($\star$).

\begin{figure}
	\includegraphics[width=1\linewidth]{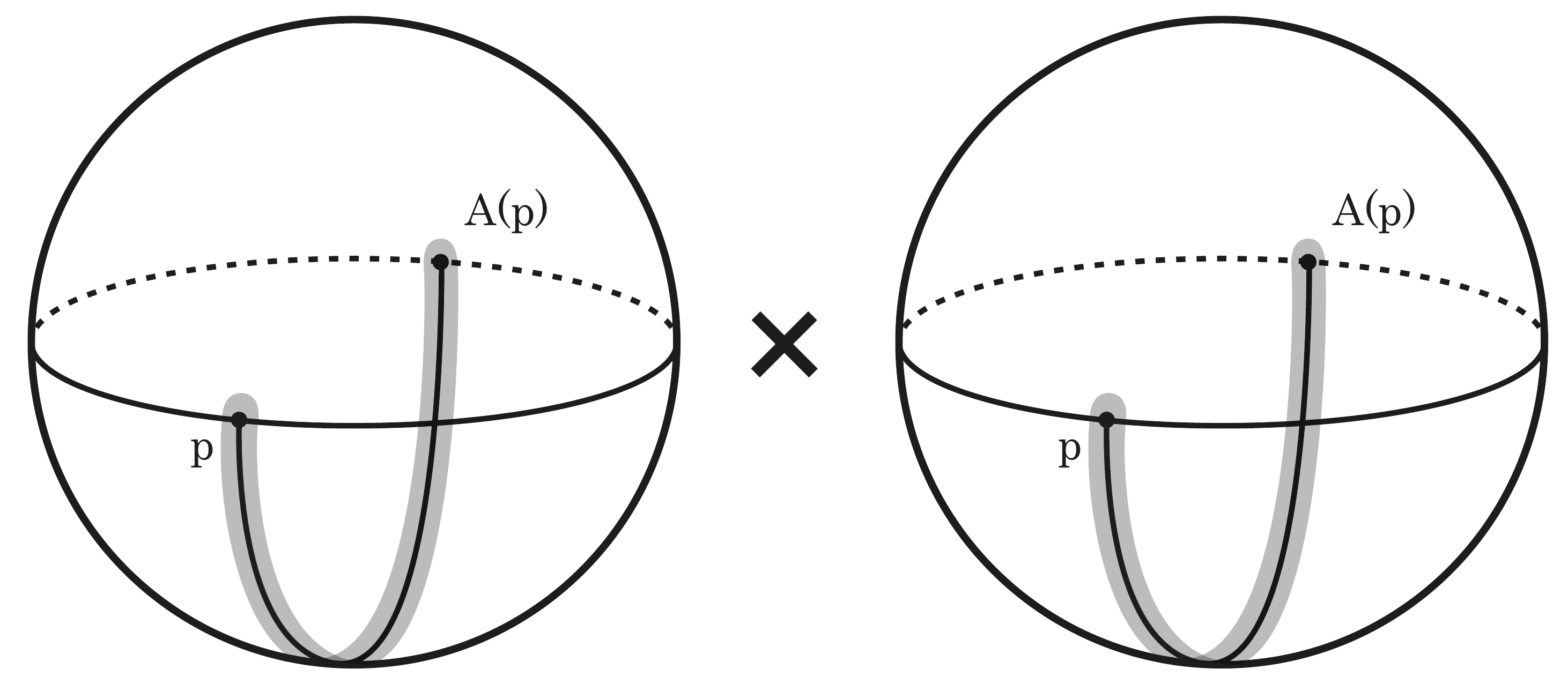}
	\caption{$(\frac{1}{2}Eq)_{\epsilon}\times (\frac{1}{2}Eq)_{\epsilon}$ remains as a barrier in $S^2\times S^2$.}
	\label{Fig4}
\end{figure}

Let us start with the simplest case, the manifold $(S^2\times S^2, \mathring{g}\times\mathring{g})$ which is isometric to $G^{+}_{2,4}$ with its standard metric. By Example~\ref{S2 without half equator}, we know that $S^{2}\backslash(\frac{1}{2}Eq)_{\epsilon}$ has property ($\star$). We define $\mathcal{R}:= \left(S^{2}\backslash(\frac{1}{2}Eq)_{\epsilon}\times (S^{2}\backslash(\frac{1}{2}Eq)_{\epsilon})\right)$ and look at a given geodesic in $S^{2}\times S^{2}$:
\begin{align}
\begin{split}
\gamma: &\R\longrightarrow S^2\times S^2\\
       &t\mapsto (\gamma^1(t),\gamma^2(t)).
\end{split}
\end{align}
where $\gamma^1,\hspace{0.1cm}\gamma^2:\R\longrightarrow S^2$ are two geodesics in $S^2$. Therefore $\mathcal{R}$ satisfies Conditions~\ref{i} and \ref{ii} of Theorem~\ref{property star} (actually, Theorem  \ref{corollary of SMP} already suffices to see that $\mathcal{R}$ is a barrier) and therefore $\mathcal{R}$ has property ($\star$) as well. (See Figure~\ref{Fig4})

By Equation~\eqref{the slope of the gauss map}, if $\Sigma^2=(x,f(x))$ is an embedded minimal surface in $\R^4$ where $f:\R^2\longrightarrow \R^2$ is a smooth mapping with slope $\Delta_{f}<\beta_0<+\infty$, then the harmonic Gauss map $\gamma:\Sigma^2\longrightarrow G^{+}_{2,4}$ has its image lying in the above region, i.e., $\gamma(\Sigma^2)\subset \mathcal{R}$. Now, by Theorem~\ref{property star for graphs}, we have that $\gamma$ is constant, implying that $\Sigma^2$ is a plane in $\R^4$.

For the more general case where $n=p+2$, we proceed as follows. Let $w\in G^{+}_{p,p+2}\subset S^{\binom{p+2}{p}-1}$ be written as $w=e_1\wedge...\wedge e_p$, where $\{e_i\}_{i=1}^{p}$ is an orthonormal basis for $V_w$ and $\{n_1,n_2\}$ an orthonormal basis for $V_{w}^{\perp}$. In $T_{w}G^{+}_{p,p+2}$ there are $\frac{p(p-1)}{2}$ vectors of type 2 (as $\wedge$ anti-commutes), we denote them by $X_{2}^{i}$, $i=1,...,\frac{p(p-1)}{2}$.

For a given $\epsilon>0$, define 
\begin{equation}
\mathcal{R}:=\displaystyle\bigcup_{i=1}^{\frac{p(p-1)}{2}}\displaystyle\big(\displaystyle\bigcup_{-t_{X^{i}_{2}}+\epsilon < t < t_{X^{i}_{2}}-\epsilon}B_{G}(w_{X^{i}_{2}}(t))\big).
\end{equation}

\noindent Since $G^{+}_{p,p+2}$ is a homogeneous space of non-negative curvature and the geodesics are given by Equation~\eqref{geodesic eq}, we have that $\mathcal{R}$ satisfies Condition~\ref{ii} in section 3, that is, for a given $\tilde{w}\in\mathcal{R}$, up to time $t<t_{X_j}$, $j=1,2$, there is always a direction where the distance to $\tilde{w}$ increases. On the other hand, as we take the union over $-t_{X_2}+\epsilon< t <t_{X_2}-\epsilon$ (i.e., strictly smaller than the critical time), we are excluding the points where $\langle w,w_{X^{i}_{2}}(t)\rangle=\cos^2(t/\sqrt{2})=0$, therefore $\langle w,w_{X_{2}^{i}}(t)\rangle>0$ for all $i=1,...,p(p-1)/2$ and there are no closed geodesics in $\mathcal{R}$ centered at any point $\tilde{w}\in\mathcal{R}$. Once more, as $G^{+}_{p,p+2}\subset S^{\binom{p+2}{p}-1}$ is a symmetric space, $w_{-X_{2}^{i}}(t)=w_{X^{i}_{2}}(-t)$, and we have that the only points in $w_{X^{i}_{2}}(\R)$ we are excluding are the ones that are $\epsilon$-close to the equator with respect to $w$ in $S^{\binom{p+2}{p}-1}$, thus $\mathcal{R}=\textstyle\big(G^{+}_{p,p+2}\backslash Hem(-w)\big)_{\epsilon}$ has property ($\star$) and hence, by Theorem~\ref{property star for graphs}, property ($\star$) for graphs, for every $\epsilon>0$.

Once more, by Equation~\eqref{the slope of the gauss map}, if $M^p=(x,f(x))$ is an embedded minimal submanifold in $\R^{p+2}$ and  $f$  is a smooth mapping with slope $\Delta_f<\beta_0 < +\infty$, then the harmonic Gauss map $\gamma: M^p\longrightarrow G^+_{p,p+2}$ satisfies $\gamma(M^p)\subset\mathcal{R}$. Thus $\gamma$ is constant and $M^p$ is a plane in $\R^{p+2}$.

\begin{rmk}
	If $n\geq p+3$, the definition of such a region $\mathcal{R}$ would have to include not only the union over directions of type 3 (and more), but unions over type 2 vectors in directions of type 1 and so on. For example, once we transport a type 1 vector along a geodesic with tangent vector of type 3, we may obtain a type 2 vector and hence we have to control with different arguments whether this would give us a closed geodesic in $\mathcal{R}$. It can be proven that in general a closed geodesic exists inside $\mathcal{R}$ if we take the union over all directions in higher codimensions. Therefore a more refined  definition of such a region in higher codimension is necessary, as follows from the existence of the  Lawson-Osserman cone.
\end{rmk}
 
\noindent Returning to codimension 2, the above discussion yields the following theorem, which extends  Moser's result from codimension 1 to codimension 2.

\begin{thm}
Let $z^i=f^i(x^1,...,x^p)$, $i=1,2$ be smooth functions defined everywhere in $\R^p$. Suppose their graph $M=(x,f(x))$ is a submanifold with parallel mean curvature in $\R^{p+2}$. Suppose that there exists a number $\beta_0 < +\infty$ such that
\begin{equation}
\Delta_{f}(x)\leq \beta_0 \hspace{0.3cm}\text{for}\hspace{0.2cm}\text{all}\hspace{0.2cm} x\in\R^n,
\end{equation}
where
\begin{equation}
\Delta_{f}(x):=\big\{\det\big(\delta_{\alpha\beta} + \sum_{i}f^{i}_{x_{\alpha}}(x)f^{i}_{x_{\beta}}(x)\big)\big\}^{\frac{1}{2}}.
\end{equation}
Then $f^1,f^2$ are linear functions on $\R^p$ representing an affine $p$-plane in $\R^{p+2}$.
\end{thm}

\bibliographystyle{alpha}
\bibliography{bernstein}
\addcontentsline{toc}{section}{\bibname}

\end{document}